\documentstyle[12pt]{amsart}

\chardef\bslash=`\\ % p. 424, TeXbook

\newtheorem{thm}{Theorem}[section]
\newtheorem{cor}[thm]{Corollary}
\newtheorem{lem}[thm]{Lemma}
\newtheorem{prop}[thm]{Proposition}

\theoremstyle{definition}
\newtheorem{defn}{Definition}[section]

\theoremstyle{remark}

\numberwithin{equation}{section}

\newcommand{\thmref}[1]{Theorem~\ref{#1}}
\newcommand{\propref}[1]{Proposition~\ref{#1}}
\newcommand{\secref}[1]{\S\ref{#1}}
\newcommand{\lemref}[1]{Lemma~\ref{#1}}
\newcommand{\corref}[1]{Corollary~\ref{#1}}
\newcommand{\defref}[1]{Definiton~\ref{#1}}

\newcommand{\WPtwo}{{{\bold P}(2,1^4)}}
\newcommand{\WP}{{\tilde{\bold P}}}

\newcommand{\oo}{{\cal O}}

\newcommand{\V}{{\cal V}}
\newcommand{\W}{{\cal W}}

\newcommand{\I}{{\cal I}}

\newcommand{\C}{{\cal C}}

\newcommand{\Hilb}{\bold{Hilb}}

\newcommand{\HilbPfour}{\bold{Hilb}_{\bold P^4}}
\newcommand{\HilbPthree}{\bold{Hilb}_{\bold P^3}}
\newcommand{\HilbWP}{\bold{Hilb}_\WP}
\newcommand{\HP}{{\cal H}}
\newcommand{\HWP}{{\tilde{\cal H}}}
\newcommand{\Spec}{\operatorname{Spec}}
\newcommand{\Proj}{\operatorname{Proj}}
\newcommand{\Grass}{\operatorname{Grass}}
\newcommand{\Coh}{\operatorname{H}}
\newcommand{\Hom}{\operatorname{Hom}}
\newcommand{\GL}{\operatorname{GL}}
\newcommand{\PGL}{\operatorname{PGL}}
\renewcommand{\phi}{\varphi}

\newcommand{\Bigdownsurjarrow}{\big\downarrow\kern-8.6pt\downarrow}

\begin{document}

\title[Rational quartics on Calabi-Yau hypersurfaces in $\WPtwo$]
{The Number of Rational Quartics on Calabi-Yau hypersurfaces in
Weighted Projective Space $\WPtwo$}

\author{Paul Meurer}
\address{Department of Mathematics\\
        University of Bergen\\ N-5007 Bergen, Norway}
\email{Paul.Meurer@@mi.uib.no}

\date{September 1, 1994}



\maketitle

\begin{abstract}
We compute the number of rational quartics on a general Calabi-Yau
hypersurface in weighted projective space $\WPtwo$. The result agrees
with the prediction made by mirror symmetry.
\end{abstract}

\section{Introduction}
\label{intro}
In this note we will compute the number of rational quartics
(see \secref{quartics})
on a general Calabi-Yau hypersurface in weighted projective space $\WPtwo$.
The number was asked for for the first time by Sheldon Katz \cite{katz:ratcur},
and David Morrison \cite{morr:picfux} computed it using mirror symmetry
methods.
Our result, which agrees with the mirror symmetry computation, is:
\begin{thm}
There are 6028452 rational quartics on a general Calabi-Yau hypersurface
in weighted projective space $\WPtwo$.
\end{thm}
The method used can be sketched as follows:
We show that the irreducible component $\HWP_4$ of the Hilbert scheme of
$\WPtwo$
containing the rational quartics is smooth and can be embedded into the
irreducible  component $\HP_4$ of the Hilbert scheme of $\bold P^4$
containing the elliptic quartic curves.
There, it can be characterized as one component of the fixed point
scheme of a natural involution.
On the other hand, $\HP_4$ is well-known and explicitly described
by Dan Avritzer and Israel Vainsencher~\cite{avvai}.
Together, this leads to an explicit description of $\HWP_4$.

The number of rational quartics on a general Calabi-Yau hypersurface is given
as
the integral of the top Chern class of a certain vector bundle on
$\HWP_4$. It will be computed by a formula of Bott's, which
expresses the integral of a homogeneous polynomial in the Chern classes
of a bundle on a smooth, compact variety with a $\bold C^*$-action in terms of
data given by the induced linear actions on the fibers of the bundle and
the tangent bundle in the (isolated) fixed points of the action.

{\it Acknowledgement.}
I would like to express my thanks to Stein Arild Str{\o}mme for many
helpful conversations.

\section{Relations between the Hilbert schemes of $X$ and $X/\Gamma$}
In this section, we investigate relations between the Hilbert schemes of
a projective scheme $X$ and its quotient $X/\Gamma = Y$ w.r.t.\ the action
of a finite group $\Gamma$. To each irreducible component $\HP$ of
$\Hilb_{X/\Gamma}$, we will find a subscheme $Z$ of
$\Hilb_X$ mapping birationally onto $\HP$. The morphism $\phi: Z\to \HP$
is generally no isomorphism, but in the situation interesting us primarily
(i.e.\ $\HP =$ irreducible component of $\Hilb_{\WPtwo}$ containing the
rational quartics), it actually is an isomorphism.\\

We assume all schemes to be defined over an algebraically closed field
$k$ of characteristic 0.

First, we state some nice simple properties of finite quotients
$\rho: X\to X/\Gamma=Y$ for later reference.
\begin{lem}\label{quotients}
\begin{enumerate}
\item $\rho_*$ and $(\rho_*(\cdot ))^\Gamma$ are exact functors
from $\Gamma$-linearized quasicoherent $\oo_X$-modules to quasicoherent
$\oo_Y$-modules.
\item Let $V\subseteq X$ be $\Gamma$-invariant. Then we have
$$\oo_{V/\Gamma} = (\rho_*\oo_V)^\Gamma,\qquad
\I_{V/\Gamma} = (\rho_*\I_V)^\Gamma.$$
\item $(\rho_*(\cdot ))^\Gamma$ commutes with cohomology, i.e.
$$\Coh^i(X,\cal F)^\Gamma = \Coh^i(Y,(\rho_*\cal F)^\Gamma).$$
\end{enumerate}
\end{lem}
\begin{pf}
(i) $\rho$ is a finite map.\\
(ii) follows immediately from the definition of $V/\Gamma$ and (i).\\
(iii) Note that the functors $\Coh^i(X,\cdot )^\Gamma$ and
$\Coh^i(Y,(\rho_*(\cdot ))^\Gamma)$ are left-exact and equal for $i=0$.
But the category of quasicoherent $\oo_X$-modules has enough injectives,
so they are equal for all $i$.
\end{pf}
\begin{lem}\label{flatquotient}
Let $X$ be a quasiprojective scheme and $\Gamma$ a finite group acting
on $X$.
Let
$$
\begin{array}{ccc}
V& \hookrightarrow& S\times X\\
\downarrow& \swarrow& \\
S& &
\end{array}
$$
be a family of subschemes of $X$, flat over $S$ and invariant under the
action of $\Gamma$.
Then $V/\Gamma \subseteq S\times X/\Gamma$ is flat over $S$.
\end{lem}
\begin{pf}
We can assume that $S$ and $X$ are affine, since flatness is a local
property and $X\to Y$ is an affine map. Write $S = \Spec A$,
$X = \Spec R_0$, $R=A\otimes R_0$.
On every $R$-module $M$, we define an $R^\Gamma$-linear endomorphism
$\Phi_M: M\to M$ by
\begin{equation}\label{average}
\Phi_M(x) := |\Gamma|^{-1}\sum_{\gamma\in\Gamma}\gamma(x).
\end{equation}
Then $M^\Gamma=\Phi_M(M)$, and $M$ splits as
$M = M^\Gamma\oplus \ker\Phi_M$.
In particular, we have $\oo_V = \tilde M$ for some $A$-flat $R$-module
$M$, so $\oo_{V/\Gamma} = \tilde{M^\Gamma}$, and $M^\Gamma$ is $A$-flat
as a direct summand of a flat $A$-module.
\end{pf}
For a finite group $\Gamma$ acting on a quasiprojective scheme $H$, we can
define the {\it fixed point scheme} $H^\Gamma$ of the $\Gamma$-action
as follows:

If $H=\Spec R$ is affine, then $H^\Gamma$ is defined by the ideal generated
by $\ker(\Phi_R)$. In general, cover $H$ by affine invariant open sets
$U$ and let $(H^\Gamma)_{|U} := U^\Gamma$.

$H^\Gamma$ can be characterized as the maximal subscheme of $H$ having
$\Gamma$-invariant structure sheaf.

Note that the group $\Gamma$ acts in a natural way on $\Hilb_X$: if
$\lbrack C\rbrack\in\Hilb_X$ and $\gamma\in\Gamma$, the action is
given by $\gamma(\lbrack C\rbrack) := \lbrack\gamma(C)\rbrack$.\\

Let now $\HP_Y$ be an irreducible component of $\Hilb_Y$ with universal
family $\W$. For any open set $U\subseteq\HP_Y$, denote by $\W_U$
the restriction of $\W$ to $U$ and let
$\rho_U := \rho\times id_U$. The lift of the family $\W_U$ to $X$ is given by
$\rho_U^{-1}\W_U = \W_U\times_{Y\times U}X\times U$.

\begin{thm}\label{birational}
Suppose that there exists an open subset $U\subseteq\HP_Y$ such that
$\rho_U^{-1}\W_U$ is flat. Then there is a uniquely determined irreducible
component $Z$ of $(\Hilb_X)^\Gamma$ mapping birationally to~$\HP_Y$ by the
map $\lbrack C\rbrack\mapsto \lbrack C/\Gamma\rbrack$.
\end{thm}
\noindent{\it Remark.} The existence of such an open set is guaranteed
for example if $\HP_Y$ is reduced (cf. \cite{mumford:curves}, Lect.~8:
Flattening stratifications).
\begin{pf}
$\rho_{\HP_Y}^{-1}\W$ is in general not flat over $\HP_Y$, but by
assumption, there is an open subset $U\subseteq \HP_Y$ such that
$\rho_U^{-1}\W_U$ is flat. We can assume that $U = \Spec A$.
By the universal property of the Hilbert functor, we
thus get a morphism $\psi: U\to \HP_X$, where $\HP_X$ is some
irreducible component of $\Hilb_X$. Moreover, $\psi$ factors through an
irreducible component $Z$ of $(\HP_X)^\Gamma$.
Let $\V_Z\subseteq X\times Z$ be the restriction of the universal family
on $\HP_X$ to $Z$.

If we denote $id_X\times\psi$ by $\psi_X$ etc., it is immediate that
$\rho_{U*}\psi_X^*\cal F=\psi_Y^*\rho_{Z*}\cal F$ for any coherent sheaf
$\cal F$ on $X\times Z$. On the other hand,
$\oo_{\W_U} = (\rho_{U*}\rho_U^*\oo_{\W_U})^\Gamma$; a simple computation shows
then that $\W_U$ is the pullback of $\V_Z/\Gamma$ by $\psi$.

 By \lemref{flatquotient}, $\V_Z/\Gamma\subseteq Y\times Z$ is flat over $Z$,
therefore we get a morphism $\phi$ from $Z$ to the same component $\HP_Y$
such that $\V_Z/\Gamma$ is the pullback of the universal family $\W_Y$.

Since $\rho_U^{-1}\W_U$ is flat over $U$,  it is again easy to see that
$\V_{\phi^{-1}U}$ is the pullback of $\rho_U^{-1}\W_U$ by $\phi$.

We can summarize the situation in the following commutative diagram:
$$
\begin{array}{ccccccccc}
\rho_U^{-1}\W_U & \longrightarrow&  \V_Z & \supseteq & \V_{\phi^{-1}U} & &
\longrightarrow & &  \rho_U^{-1}\W_U\\
\downarrow&  &\downarrow& & & & & & \downarrow\\
\W_U&\longrightarrow &\V_Z/\Gamma& &\longrightarrow & & \W & \supseteq & \W_U\\
\downarrow&  &\downarrow& & & &\downarrow& &\downarrow\\
U &\stackrel{\psi}{\longrightarrow} & Z & &\stackrel{\phi}{\longrightarrow} & &
\HP_Y & \supseteq & U
\end{array}
$$
Since all maps in the first and second rows are pullbacks by $\psi$
or $\phi$, it follows by the universal property of the Hilbert functor
that $\phi\circ\psi=id_U$ and $\psi\circ\phi_{\phi^{-1}U}=id_{\phi^{-1}U}$.
Therefore, $\phi: Z\to\HP_Y$ is birational.
\end{pf}
{\it Example.} A simple example showing that the map $\phi:Z\to\HP_Y$ dosen't
need to be an isomorphism is the following:

Take $X=\bold P^2$ and $\Gamma = \{id,\iota\}$, where $\Gamma$ acts on
$\bold P^2$ by $\iota x_0=-x_0$, and $\iota x_i=x_i$ for $i=1,2$.
Then $Y=X/\Gamma$ is the weighted projective space $\bold P(2,1,1)$
(cf. \defref{def:weighted_pr_sp}).
Let $\HP^1_{\bold P(2,1,1)}$ be the Hilbert scheme of subschemes of length 1 in
$\bold P(2,1,1)$, which is isomorphic to $\bold P(2,1,1)$ itself, and singular
in the point $(1,0,0)$. On the other hand, since the inverse image of a general
point
in $\bold P(2,1,1)$ is a pair of points in $\bold P^2$, $Z$ lies in the Hilbert
scheme $\HP^2_{\bold P^2}$ of subschemes of length 2 in $\bold P^2$, which is a
$\bold P^2$-bundle over $\check{\bold P}^2$, hence smooth. By
\corref{cor:smooth},
$Z$ is smooth too, hence $Z$ can't be isomorphic to $\bold P(2,1,1)$.
In fact, $Z$ is the blow up of  $\bold P(2,1,1)$ in the singular point.\\

The following proposition is certainly well known. Since we didn't find
a reference for it, we give a proof.
\begin{prop}
Let $H$ be a smooth scheme and $\Gamma$ a finite group acting on $H$.
\begin{enumerate}
\item The fixed point scheme $H^\Gamma$ of $\Gamma$ is a smooth
subscheme of $H$.
\item The Zariski tangent space $\cal T_{H^\Gamma}(x)$ to $H^\Gamma$
in a point $x$ is equal to $\cal T_H(x)^\Gamma$.
\end{enumerate}
\end{prop}
\begin{pf}
(i)
Let $p\in H^{\Gamma}$ be a closed point and $R = \oo_{p,X}$; consider the
induced action of $\Gamma$ on $R$.
Then the ideal $\bold a_\Gamma$ of $H^{\Gamma}$ in $R$ is generated by
$\ker\Phi_R$.

Since $H$ is smooth, $R$ is a regular (noetherian) local ring. We will
show that also $R/\bold a_\Gamma$ is regular, which proves the proposition.

Let $\bold m\subseteq R$ be the maximal ideal. $\Gamma$ induces a linear
action on the vector space $\bold m/\bold m^2 = V$.

Let $\hat x_1,\ldots,\hat x_d$ be a base of the subspace of $V$ invariant under
$\Gamma$, and complete it by $\hat y_1,\ldots,\hat y_{n-d}$ to a
base of $V$, such that the $\hat y_j$ satisfy $\Phi_R(\hat y_j)=0$.
This can be achieved by choosing arbitrary $y_j$ and setting
$\hat y_j=y_j-\Phi_R(y_j)$. By Nakayama's lemma, the $\hat x_i$,
$\hat y_j$ lift to generators $\tilde x_i$, $\tilde y_j$ of $\bold m$.

Again we can construct new ring elements by averaging:
Let $x_i := \Phi_R(\tilde x_i)$, $y_j := \tilde y_j-\Phi_R(\tilde y_j)$.
The images of $x_i$, $y_j$ in $R/\bold m^2$ are $\hat x_i$, $\hat y_j$.
We conclude that also $x_i$, $y_j$ generate $\bold m$.
We will show that $\bold a := (y_1,\ldots,y_{n-d})$ is equal to
$\bold a_\Gamma$; then $R/\bold a_\Gamma$ is a regular local ring,
and $H^\Gamma$ is smooth.

The ideal $\bold a$ is clearly contained in $\bold a_\Gamma = (\ker\Phi_R)$.
On the other hand, suppose that
$$\ker\Phi_R\subseteq \bold a + \bold m^r,$$ (this is trivially true for
$r=1$), and let
$b=\sum_i\alpha_ix_i+\sum_j\beta_jy_j \in\ker\Phi_R$ be an arbitrary
element. We can write
$$b = b-\Phi_R(b) =
\sum (\alpha_i-\Phi_R(\alpha_i))\,x_i + \sum\beta_jy_j-\Phi_R(\beta_jy_j).$$
The first sum is contained in $\bold a+\bold m^{r+1}$ because
$\alpha_i-\Phi_R(\alpha_i)$ lies in $\ker\Phi_R\subseteq\bold a +\bold m^r$.
The summands $\beta_jy_j$ lie in $\bold a$.
Write now
\begin{eqnarray*}
\Phi_R(\beta_jy_j) &=
&\Phi_R\big(\Phi_R(\beta_j)\,y_j+(\beta_j-\Phi_R(\beta_j))\,y_j\big)\\
&=&\Phi_R((\beta_j-\Phi_R(\beta_j))\,y_j)\\
&=&|\Gamma|^{-1}\sum_{\gamma}\gamma(\beta_j-\Phi_R(\beta_j))\,\gamma(y_j).
\end{eqnarray*}
In the last sum, both factors of the summands lie in $\ker\Phi_R$, hence
$\Phi_R(\beta_jy_j)$ lies in $\bold a +\bold m^{r+1}$, too, and together we
have
that $b\in\bold a +\bold m^{r+1}$, hence
$$\ker\Phi_R\subseteq \bold a + \bold m^{r+1}.$$
Since $\bigcap_r\bold m^r=\emptyset$, we get by induction on $r$ that
$$\ker\Phi_R\subseteq\bold a.$$
Hence $\bold a_\Gamma = \bold a$, and $H^\Gamma$ is smooth.

(ii) is immediate.
\end{pf}

With the notation of \thmref{birational}, we have
\begin{cor}\label{cor:smooth}
$Z$ is smooth if $\HP_X$ is smooth.
\end{cor}

\section{Rational quartics in $\WPtwo$ and elliptic quartics in $\bold P^4$}
\label{quartics}
\begin{defn}\label{def:weighted_pr_sp}
A {\em weighted projective $n$-space} $\bold P(k_0,\ldots, k_n)$
with positive integer weights $k_0,\ldots,k_n$ is defined
as $\Proj$ of the graded ring $R := k\lbrack y_0,\ldots,y_n\rbrack$,
the variables $y_0,\ldots,y_n$ having weights $k_0,\ldots,k_n$.
We can assume that the $k_0,\ldots,k_n$ are coprimal.

There is an isomorphism of graded rings
$R\cong k\lbrack x_0^{k_0},\ldots,x_n^{k_n}\rbrack$, the $x_i$ having weight 1,
and in the following we will work with the $x_i$ rather than the $y_i$.
\end{defn}
We will use the abbreviation $\WP$ for an arbitrary weighted projective
space.

\smallskip
Let  $\Gamma_k$ be the group of $k$-th roots of unity, and let
$\Gamma := \Gamma_{k_0}\times\ldots\times\Gamma_{k_n}$.
$\Gamma$ acts on $\bold P^n= \Proj k(x_0,\ldots,x_n)$ in the obvious way,
namely through the action
of $\Gamma_{k_i}$ on the $i$-th projective coordinate by multiplication.
$\WP$ can then be described as the geometric quotient of $\bold P^n$
by this action of $\Gamma$:
$$\bold P(k_0,\ldots,k_n) = \bold P^n/\Gamma.$$
We denote by $\rho$ the quotient map $\bold P^n\to\WP$, which is a
finite ramified covering map, induced by the inclusion
$k\lbrack x_0^{k_0},\ldots,x_n^{k_n}\rbrack\subset k\lbrack
x_0,\ldots,x_n\rbrack$.
The ramification locus of $\rho$ is the union
of the fixed point sets of the non-identity elements of $\Gamma$.

The singular points of $\WP$ are exactly the points
$\stackrel{i}{(0,\ldots,0,1,0,\ldots,0)}$ such that $k_i>1$.

We can define twisting sheaves in the usual way by setting
$\oo_\WP(r) := R(r)\tilde{\ }$. In general, these sheaves need not to be
invertible, but when $\WP = \bold P(k,1^n)$, then $\oo_\WP(k)$ is a very ample
line bundle and induces a Veronese type embedding of $\WP$ into a big
projective space $\bold P^N$: the image of $\WP$ is the projective cone
over the image $v_k(\bold P^{n-1})$ of the Veronese embedding of
$\bold P^{n-1}$ by $\oo_{\bold P^{n-1}}(k)$.

There is a natural projection from the singular point onto $\bold P^{n-1}$
induced by the inclusion
$k\lbrack x_1,\ldots,x_n\rbrack\subset k\lbrack x_0^k,x_1,\ldots,x_n\rbrack$.
In $\bold P^N$, this projection is given by projecting the cone down to
$v_k(\bold P^{n-1})$ from the point $(1,0,\ldots,0)$.\\

We will need the following
\begin{lem} The quotient map $\rho:\bold P^n\to\bold P(k,1^n)$ is flat away
from the singular point $(1,0,\ldots,0)$.
\end{lem}
\begin{pf}
$\bold P^n$ can be covered by the $\Gamma$-invariant affine open sets
$D_+(x_i)$, $i=0,\ldots,n$, and the union
$\bigcup_{i=1}^n D_+(x_i)$ is equal to $\bold P^n-\{(1,0,\ldots,0)\}$.
It suffices to show that all the maps $D_+(x_i)\to D_+(x_i)/\Gamma$ are
flat for $i=1,\ldots,n$.
Let $$R_i= k\Big\lbrack{x_0\over x_i},\ldots,{x_{i-1}\over x_i},
{x_{i+1}\over x_i},\ldots,{x_n\over x_i}\Big\rbrack.$$
Then the invariant ring is
$$R_i^\Gamma= k\Big\lbrack\Big({x_0\over x_i}\Big)^{\!k},{x_1\over x_i},\ldots,
{x_{i-1}\over x_i},{x_{i+1}\over x_i},\ldots,{x_n\over x_i}\Big\rbrack,$$
and we have $D_+(x_i)=\Spec R_i$, $D_+(x_i)/\Gamma= \Spec R_i^\Gamma$.
Thus the ring $R_i$ is equal to
$\bigoplus_{s=0}^{k-1}({x_0\over x_i})^s\cdot R_i^\Gamma$,
hence a free $R_i^\Gamma$-module, and $D_+(x_i)$ is flat over
$D_+(x_i)/\Gamma$.
\end{pf}
\begin{cor}\label{flatlift}
Let $V\hookrightarrow S\times\bold P(k,1^n)$ be a flat family over $S$ such
that no fiber of $V$ contains the singular point of $V$. Then
$\rho_S^{-1}V\subseteq S\times\bold P^n$ is flat over $S$.
\end{cor}
\begin{pf}
$V$ is a flat family in $\bold P(k,1^n)-\{(1,0,\ldots,0)\}$. The claim
follows from the lemma and by transitivity and base change stability of
flatness.
\end{pf}

There is no intrinsic notion of degree of a curve in a weighted projective
space $\WP$. In the case of a weighted projective space of type
$\WP = \bold P(k,1^n)$, however, we agree to measure the degree
of curves (or, more generally, the Hilbert polynomial of subschemes)
w.r.t.\ the embedding described above, i.e.\
$\deg C =\deg(\oo_\WP(k)_{|C})$.\\

Let us now specialize to the case $\WP = \WPtwo$, the case we are primarily
interested in.
$\WP=\Proj(k\lbrack x_0^2,x_1,\ldots,x_4\rbrack)$ is the quotient of $\bold
P^4$
by the group $\Gamma=\{1,\iota\}\cong\bold Z_2$, the action of $\Gamma$
being given by $\iota x_0=-x_0$, $\iota x_i=x_i$ for $i=1,\ldots,4$.

We will consider rational quartics, i.e.\ rational curves of degree
(in the above sense) four in $\WP$. Let $C\subseteq\WP$ be a rational
quartic and
\begin{eqnarray}\label{parametrization}
	\bold P^1 & \to & C\nonumber \\
	(s,t) & \mapsto & \big(f_0(s,t),\ldots,f_4(s,t)\big)
\end{eqnarray}
a parametrization of $C$. The image of $C$ under the embedding
$\WP\hookrightarrow\bold P^N$ induced from $\oo_\WP(2)$ is parametrized as
$$(s,t)\quad\mapsto\quad (f_0,f_1^2,f_1f_2,\ldots,f_4^2)$$
and is a degree four curve by definition. Hence $f_0$ has degree four and
$f_1,\cdots,f_4$ have degree two. Furthermore, we see that the projection
of a rational quartic from the singular point to $\bold P^3$ is a conic
in $\bold P^3$.

We denote by $\HWP_4$ the irreducible component of $\HilbWP$ containing
the rational quartics.
(To see that the rational quartics are really contained in {\it one}
irreducible
component, observe that the parameter space of rational quartic curves
(without degenerations) is fibered
over the space of conics in $\bold P^3$ (which is irreducible). Each fiber is
irreducible, too, as one can seefrom \eqref{parametrization} by specifying a
parametrization $(f_1,\ldots,f_4)$ of a conic and letting $f_0$ vary.
Hence the space of rational quartics is irreducible, and $\HWP_4$ is the
irreducible component of $\HilbWP$ containing that space.)
The universal curve on $\HWP_4$ will be denoted by $\C$.

Let now $C$ be a general rational quartic in $\HWP_4$.
The quotient map $\rho^{-1}C\to\rho^{-1}C/\Gamma=C$ exhibits
$\rho^{-1}C$ as a double cover of $C$, and the ramification locus
consists of the four points where $f_0=0$ in the parametrization
\eqref{parametrization}. By Hurwitz's theorem, $\rho^{-1}C$ is an
elliptic curve, and it is again a quartic because $\rho^{-1}C$ intersects
the hyperplane $\{x_0=0\}$ in four points.

Since $C$ doesn't contain the singular point of $\WP$, there is an open subset
$U\subseteq\HWP_4$ containing $\lbrack C\rbrack$ such that no fiber of $\C_U$
contains the singular point. By \corref{flatlift}, $\C_U$ lifts to a flat
family  $\rho_U^{-1}\C_U$ in $\bold P^4$.

Denote by $\HP_4$ the component of $\HilbPfour$ parametrizing the smooth
elliptic
quartic curves in $\bold P^4$ and their degenerations.
According to \thmref{birational}, there is an irreducible subscheme $Z_4$
of $\HP_4$ mapping birationally to $\HWP_4$ by
$\lbrack C'\rbrack\mapsto \lbrack\rho(C')\rbrack$. Since $\HP_4$ is smooth
(see below), it follows by \corref{cor:smooth} that $Z_4$ is smooth, too.\\

\noindent{\it Remark.} It is clear that the above considerations are valid as
well
for rational quartics in $\bold P(2,1^3)$ instead of $\bold P(2,1^4)$.
Thus, if we denote by $\HP_3$ the irreducible component of $\HilbPthree$
parametrizing the elliptic quartic curves in $\bold P^3$,
there is a smooth irreducible subscheme $Z_3$ of $\HP_3$ mapping birationally
to $\HWP_3$ by $\lbrack C'\rbrack\mapsto \lbrack\rho(C')\rbrack$.\\

In order to obtain the explicit description of $Z_4$ that we need for the
calculations, and particularly to prove the claimed isomorphy of $Z_4$ and
$\HWP_4$, we will have to look at the description of $\HP_4$ given by
Dan Avritzer and Israel Vainsen\-cher:\\

Let $G := \Grass_2(\Coh^0(\oo_{\bold P^3}(2)))$
be the Grassmannian of pencils of quadric surfaces in $\bold P^3$, and
denote by $G'$ the image of the (well-defined) map
\begin{eqnarray}\label{H3toG'}
\HP_3 &\to &\Grass_8(\Coh^0(\oo_{\bold P^3}(3))\\
\lbrack C \rbrack &\mapsto &\lbrack \Coh^0(\I_C(3))\rbrack\ .\nonumber
\end{eqnarray}
On $G$, we have a canonical family of subschemes of $\bold P^3$: the fiber
in a point $g\in G$ is the base locus of the
pencil represented by $g$. In the same way, $G'$ gives rise to a family
of subschemes of $\bold P^3$: the fiber in a point $g'$ is the base locus
of the linear system of cubic surfaces represented by $g'$.

Denote by $B$ the subscheme of $G$ where the family on $G$ is not flat,
and denote by $D$ the subscheme of $G'$ where the family defined by $G'$
is not flat ($B$ consists of pencils with a fixed component, and $D$ is
the scheme of planes in $\bold P^3$ with an embedded subscheme of length 2).
Then we have:
\goodbreak
\begin{thm}\label{thm:elliptic}
\begin{enumerate}
\item (Avritzer, Vainsencher \cite{avvai})
$\HP_3$ is isomorphic to a two-fold blow up of $G$. More precisely,
$G'$ is isomorphic to the blow up of $G$ along $B$ and $\HP_3$
is isomorphic to the blow up of $G'$ along $D$.
The ideal of every (degenerated) curve in $\HP_3$ is generated in
degree $\leq 4$. In particular, $\HP_3$ is smooth of dimension~16.
\item
$\HP_4$ is fibered locally trivially over $\check{\bold P}^4$ in a
natural way with fiber
$\HP_3$; i.e., the restriction of the universal curve over $\HP_4$
to the fiber over a hyperplane $h\cong\bold P^3$ in $\check{\bold P}^4$
is the universal elliptic quartic curve in $h\cong\bold P^3$.
\end{enumerate}
\end{thm}
\begin{pf} (ii)
We have to show that all degenerations of elliptic quartic curves in
$\bold P^4$ span exactly a $\bold P^3$. Then the fibration is given
by projecting a point $\lbrack C \rbrack\in\HP_4$ to the hyperplane it
spans.

It is clear that no degeneration can possibly span less than a $\bold P^3$,
for then this degeneration would already be contained in $\HP_3$.
On the other hand, for a general elliptic curve, $\dim(\Coh^0(\cal I_C(1)))=1$;
if a degeneration $C_0$ spanned all of
$\bold P^4$, then $\dim(\Coh^0(\cal I_C(1)))$ would drop to 0 in $C=C_0$,
in contradiction to upper semicontinuity of $\dim(\Coh^0(\cal I_C(1)))$.

The fibration is locally trivial because $\PGL(4)$ acts transitively on
$\bold P^4$ and this action lifts to an action on $\HP_4$.
\end{pf}

As a corollary of this theorem we are now able to derive an analogous
description of $Z_3$ and $Z_4$.
Consider therefore the inclusion of grassmannians
$$\tilde G:= \Grass_2\bigl(\Coh^0(\oo_{\bold P(2,1^3)}(2))\bigr)
\hookrightarrow G = \Grass_2\bigl(\Coh^0(\oo_{\bold P^3}(2))\bigr)$$
induced by the natural inclusion
\begin{equation}\label{inclusion}
i:\Coh^0(\oo_{\bold P(2,1^3)}(2))\hookrightarrow\Coh^0(\oo_{\bold P^3}(2)).
\end{equation}
%\goodbreak
\begin{prop}\label{prop:hilbertscheme}
$Z_3$ is smooth of dimension 10 and isomorphic to the proper transform of
$\tilde G$ under the twofold blow up map $b:\HP_3\to G$.
\end{prop}

\begin{pf}
$\Gamma$ acts on $G$. If $Z'$ is an irreducible component of $G^\Gamma$
not contained in the blow up locus $B$, then the proper transform of $Z'$ is an
irreducible component of $\HP_3^\Gamma$. On the other hand, it is easy
to see that $\tilde G\subseteq G$ {\it is} a component of $G^\Gamma$.
Furthermore, curves in $\tilde G-B$ map to rational quartics; this
proves that the proper transform of $\tilde G$ in $\HP_3$ is the right
component, i.e.\ equal to $Z_3$.
\end{pf}
We turn to the explicit description of $Z_3$ as a twofold blow up of
$\tilde G$.

A pencil in $\tilde B = \tilde G \cap B$ is a pencil generated by
two quadratic polynomials $F_1,F_2$ in $\Coh^0(\oo_{\WP^3}(2))$
with a common linear factor, thus we have $F_1 = f_1g$, $F_2 = f_2g$, and they
must be independent of $x_0$. It follows that the scheme described by such
a pencil projects to the union of a line and a point in
$\bold P^2\cong \{x_0=0\}$ (or a degeneration thereof) under projection
from the singular point. The image is described by the same equations
$F_1 = f_1g$, $F_2 = f_2g$. We conclude that $\tilde B$ is isomorphic to
$\check{\bold P}^2\times\bold P^2$ and has dimension 4.

According to \cite{avvai}, a plane in $\bold P^3$ with an embedded
subscheme of length 2 has ideal
$$(x_1(x_1-x_4)x_3, x_2x_3, x_3^2)
\qquad\hbox{ or }\qquad
(x_1^2x_3, x_2x_3, x_3^2)$$
up to projective equivalence, depending on whether the subscheme has
support in two points or one point.
Consequently, ideals in $\tilde D$, the intersection of $D$ with the
proper transform $\tilde G'$ of $\tilde G$,
have, up to $\Gamma$-invariant projective equivalence, the form
$$(x_3(ax_0^2+bx_1^2), x_3x_2, x_3^2)\ ,\qquad (a,b)\in \bold P^1.$$
Geometrically, $\tilde D$ parameterizes hyperplanes
through $(1,0,0,0)$ with an embedded $\Gamma$-invariant subscheme of
length 2.
By projecting the hyperplane down to $\bold P^2\cong \{x_0=0\}$,
we see that $\tilde D$ is fibered over $\check{\bold P}^2$, with fiber
$F$ isomorphic to $\bold P^1\times\bold P^1$.
If $F$ is the fiber over the point $\{x_3=0\}$, the isomorphism is given by
\begin{eqnarray*}
\bold P^1\times\bold P^1\!\quad &\to &\quad F\\
(a,b)\times(c,d) &\mapsto &(x_3(ax_0^2-b(cx_1+dx_2)^2),
x_3(dx_1-cx_2),x_3^2)\ .
\end{eqnarray*}
So $\tilde D$, too, is smooth of dimension 4.\\

\noindent{\it Remark.} This is another proof of the smoothness of $Z_3$.

\begin{prop}\label{prop:isomorphy}
The map $\phi:Z_3\to\HP_3$ is an isomorphism.
\end{prop}

\begin{pf}
Since $\phi$ is birational and $\HWP_3$ is irreducible, it is enough to
show that the tangent map $d\phi:\cal T_{Z_3}([C])\to\cal
T_{\HWP_3}([C/\Gamma])$
is injective in every point.

Consider a $\bold C^*$-action on $\bold P^3$ acting diagonally w.r.t.\
the coordinates $x_0,\ldots,x_3$ and having isolated fixed points. We can
choose the action in such a way that the induced action on $\HP_3$ has
isolated fixed points, the fixed points being given by monomial ideals.
Furthermore, the action leaves $Z_3$ invariant, hence descends to an
action on $\HWP_3$ (more about torus actions in \secref{calc}).

Suppose now that $[C]\in Z_3$ is a point where $d\phi$ is not injective.
Then $d\phi$ is not injective in any point of $\bold C^*\!\cdot\! [C]$.
Let $[C_0]\in\overline{\bold C^*\!\cdot\! [C]}$ be a fixed point in the
closure.
$d\phi$ cannot be injective in $[C_0]$ neither, because then it would be
injective on a whole neighborhood of $[C_0]$ which would contain
also points of $\bold C^*\!\cdot\! [C]$.

Therefore it suffices to show that  $d\phi$ is injective in all points
of $Z_3$ defined by monomial ideals.
A sufficient condition for this to happen is given by the following
\begin{lem}
Let $\lbrack C \rbrack\in Z_3$ be a point whose homogeneous ideal $I_C$
is generated by invariant monomials. Assume that there is no invariant
monomial $m$ of degree 3 not contained in $I_C$ such that all
monomials $x_1m,x_2m,x_3m$ lie in $I_C$.
Then $d\phi$ is injective in $\lbrack C \rbrack$.
\end{lem}
We can check explicitly by looking at all types of monomial
ideals, which we will determine in \secref{calc}, that the assumption
of the lemma is always met. This is a boring but simple exercise,
and we will only point out the reasoning for one case.

Take for example the fixed curve $C$ with ideal
$I_C= (x_1x_2,x_1x_2,x_0^2x_2)$, and let $m$ be a monomial of degree
3 not contained in $I_C$. Since $m$ is invariant, there are two
possibilities:
\begin{enumerate}
\item $x_0^2$ divides $m$. Then we have $m=x_0^2x_i$, $i=1$ or $3$, and
$x_im=x_0^2x_i^2$ is certainly not contained in $I_C$.
\item $x_0^2$ does not divide $m$. We have either $x_1\not\:\mid m$, then
$x_2m$ or $x_3m \not\in I_C$, because not both can be multiples of
$x_0^2x_2$; or $x_1\mid m$, then $m=x_1^3$, and $x_1m\not\in I_C$.
\end{enumerate}
All other cases can be treated in a similar way. Thus $d\phi$ is
injective everywhere, and the proposition is proved.

It remains to prove the lemma.

Let $[C]\in Z_3$ be a point with homogeneous ideal
$I_C=(h_1,\ldots,h_r)$, the $h_i$ being $\Gamma$-invariant monomials
of degree $d_i$ (i.e., they contain $x_0$ in even power).

There is a natural injection
\begin{equation}
\cal T_{\HP_3}([C])\hookrightarrow\Hom(\I_C, \oo_C).
\end{equation}
The resolution of $\I_C$ defined by $I_C$ gives rise to a diagram
\begin{equation}
\begin{array}{ccccc}
\bigoplus_i\oo_{\bold P^3}(-d_i)& \stackrel{I_C}{\longrightarrow}& \I_C
&\rightarrow & 0\\
&\searrow\bar f &\downarrow f & & \\
& &\oo_C & &
\end{array}
\end{equation}
for every $f\in\Hom(\I_C,\oo_C)$, and hence to an injective map
\begin{equation}
\Hom(\I_C,\oo_C) \stackrel{j}{\hookrightarrow}
\bigoplus_i\Hom(\oo_{\bold P^3}(-d_i),\oo_C)\cong
\bigoplus_i\Coh^0(\oo_C(d_i))
\end{equation}
which sends $f$ to $\bar f$.
Now let $I'_C = (h'_1,\ldots, h'_{r'})$ be a homogeneous ideal of $C$
generated by monomials $h'_i$ of {\it even} degree. In concrete terms,
construct $I'_C$ from $I_C$ by retaining the generators of even degree
and replacing each generator $h_i$ of odd degree by generators
$x_0h_i,\ldots, x_3h_i$. It is clear that $I_C$ and $I'_C$ both define
the same scheme $C$.
Let
\begin{equation}\label{I_C}
\bigoplus_{i,\lambda}\oo_{\bold P^3}(-d'_{i,\lambda})
\stackrel{I'_C}{\to}\I_C\to 0
\end{equation}
be the corresponding surjection, where $d'_{i,\lambda} = d_i+1$ and
$\lambda = 0,\ldots, 3$ if $d_i$ is odd and $d'_{i,\lambda} = d_i$
and $\lambda = -1$ (say) if $d_i$ is even.
By applying the exact functor $(\rho_{X*}(\cdot ))^\Gamma$ to \eqref{I_C},
we get a surjective map
$$\bigoplus_{i,\lambda}\oo_\WP(-d'_{i,\lambda})
\stackrel{(\Phi(h'_{i,\lambda}))_{i,\lambda}}{\longrightarrow}
\I_{C/\Gamma}\to 0.
$$
But when $d_i$ is odd, then $\Phi(h'_{i,0})=\Phi(x_0h_i)=0$, and we
don't lose anything if we let the sum run only over indices $(i,\lambda)$
with $\lambda\ne 0$.

The diagram
\begin{equation}
\begin{array}{ccccc}
\bigoplus_{i,\lambda}\oo_{\bold P^3}(-d'_{i,\lambda}) &
\stackrel{I'_C}{\longrightarrow} &\I_C & \to & 0\\
\downarrow\mu & & \Arrowvert& & \\
\bigoplus_i\oo_{\bold P^3}(-d_i)& \stackrel{I_C}{\longrightarrow}& \I_C
&\rightarrow & 0
\end{array}
\end{equation}
induces a diagram
\begin{equation}
\begin{array}{ccc}
\Hom(\I_C,\oo_C) & \stackrel{j}{\hookrightarrow}
&\bigoplus_i\Coh^0(\oo_C(d_i)) \\
\Arrowvert & & \downarrow \kappa \\
\Hom(\I_C,\oo_C) & \stackrel{j'}{\hookrightarrow}
&\bigoplus_{i,\lambda}\Coh^0(\oo_C(d'_{i,\lambda})).
\end{array}
\end{equation}
The map $\kappa$ sends an element $\ (0,\ldots,\stackrel{i}{m},\ldots, 0)\ $
to $\ (0,\ldots,\stackrel{(i,-1)}{m},\ldots, 0)\ $ if $d_i$ is even and to
$\ (0,\ldots,\stackrel{(i,0)}{x_0m},\ldots,\stackrel{(i,3)}{x_3m},
\ldots, 0)\ $
otherwise.

The tangent space to $Z_3$ in $[C]$ is equal to $(\cal T_{\HP_3}([C])^\Gamma$,
and there is a commutative diagram
\begin{equation}
\begin{array}{ccc}
\cal T_Z([C]) & \hookrightarrow &\Hom(\I_C,\oo_C)^\Gamma \\
\downarrow d\phi & & \downarrow \sigma \\
\cal T_{\HWP_3}([C/\Gamma]) & \hookrightarrow &
\Hom((\rho_*\I_C)^\Gamma,(\rho_*\oo_C)^\Gamma)
\end{array}
\end{equation}
(recall that $\I_{C/\Gamma} = (\rho_*\I_C)^\Gamma$ etc.).
$\sigma$ is induced from the functor $(\rho_{X*}(\cdot ))^\Gamma$.

The last two diagrams fit together into a big diagram
\begin{equation}
\begin{array}{ccccc}
\cal T_Z([C]) & \hookrightarrow &\Hom(\I_C,\oo_C)^\Gamma &
\stackrel{j}{\hookrightarrow}&\bigoplus_i\Coh^0(\oo_C(d_i))^\Gamma  \\
&&& \searrow & \downarrow\kappa\\
\big\downarrow d\phi &&\big\downarrow\sigma &&
\bigoplus_{i,\lambda}\Coh^0(\oo_C(d'_{i,\lambda}))^\Gamma\\
&&&&\downarrow\tau\\
\cal T_{\HWP_3}([C/\Gamma]) & \hookrightarrow &
\Hom((\rho_*\I_C)^\Gamma,(\rho_*\oo_C)^\Gamma)&\hookrightarrow &
\bigoplus_{i,\lambda\atop\lambda\ne 0}
\Coh^0(\oo_C^\Gamma(d'_{i,\lambda})).\\
\end{array}
\end{equation}
The map $\tau$ simply forgets the components with indices $(i,0)$
and is an isomorphism on the complement.

Now, in order to prove that $d\phi$ is injective, it suffices to show
that the composition $\tau\circ\kappa$ is injective.

By computing resolutions of all the ideals $\I_C$ (in Macaulay,
for example), we see that they all are at least 3-regular,
which means that
$\Coh^p(\I_C(m-p)) = 0$ for all $p>0$, $m\geq 3$; thus
$\Coh^1(\I_C(d))=0$ for $d\geq 2$. Therefore
$\Coh^0(\oo_C(d_i))^\Gamma$ is generated by the invariant monomials of
$\Coh^0(\oo_{\bold P^3}(d_i))$ not contained in
$\Coh^0(\I_C(d_i))^\Gamma$.

Thus it is enough to check that no such monomial is contained
in $\ker(\tau\circ\kappa)$. For monomials of even degree this is
clear, because $\tau\circ\kappa$ is the identity on the direct summands
$\Coh^0(\oo_C(d_i))^\Gamma$ for $d_i$ even.
An invariant monomial $m$ (which for short will stand for
$(0,\ldots,m,\ldots,0)$) of odd degree $d_i$ ($d_i=3$ in our case for all
ideals) is mapped to $[(0,\ldots,x_1m,x_2m,x_3m,\ldots,0)]$ by
$\tau\circ\kappa$. Thus $\tau\circ\kappa(m)$ is zero exactly when
$x_1m,x_2m,x_3m$ all lie in $\Coh^0(\I_C(d_i+1))^\Gamma$ (or, what
amounts to the same, in $I_C$).
\end{pf}

Now we are able to deduce the analogue to \thmref{thm:elliptic} (ii).
\begin{prop}\label{fiberingofhilb}
$Z_4$ is isomorphic to $\HWP_4$, and $\HWP_4$ is a locally trivial
fibration over $\check{\bold P}^3$ with fibers isomorphic to
$\HWP_3$. Thus $\HWP_4$ is smooth of dimension 13.
\end{prop}
\begin{pf}
We can proceed as in the proof of \thmref{thm:elliptic} (ii).

A smooth rational quartic $C$ in $\WPtwo$ spans exactly a hyperplane
$H\cong\bold P(2,1^3)$ (i.e.\ a hypersurface with a linear equation
$l(x_1,\ldots,x_4)$; use the fact that the lift of $C$ to $\bold P^4$
is an invariant elliptic quartic).
Thus, the projection of $C$ to $\bold P^3\subseteq\bold P(2,1^3)$ from
the singular point spans a hyperplane in $\bold P^3$, i.e.\ is
a point in $\check{\bold P}^3$.

Again, no degeneration can possibly span less than a hyperplane because
it then would already be contained in $\HWP_3$. But every degeneration
in $\HWP_3$ spans the whole $\bold P(2,1^3)$.
A semicontinuity argument shows again that no degeneration can span
more than a hyperplane. Thus the map
\begin{eqnarray*}
\HWP_4 &\to &\check{\bold P}^3\\
\lbrack C\rbrack &\mapsto & [span(pr(C))]
\end{eqnarray*}
is welldefined, and clearly a locally trivial fibration.

Consider now the map
$$
Z_4\to\HWP_4\to\check{\bold P}^3
$$
which exhibits $Z_4$ as a locally trivial fibration over
$\check{\bold P}^3$, with fiber $Z_3$. But $Z_3$ is isomorphic to
$\HWP_3$, so $Z_4$ is isomorphic to $\HWP_4$.
\end{pf}

\section{Rational quartics on Calabi-Yau hypersurfaces}
{}From now on, we work over the ground field $\bold C$.

A hypersurface in $\WP = \WPtwo$ given by a polynomial
of weighted degree 6 has trivial canonical bundle, i.e., is Calabi-Yau
(the following is the only point where the Calabi-Yau comes into play).

Consider the ideal sequence of the universal family $p:\C\to\HWP_4$
of rational quartics in $\WP$:
$$0\to\cal I_\C\to \oo_{\HWP_4\times\WP}
\to\oo_\C\to 0\ .$$
By twisting with $\oo(6)$ and taking direct images under $p$, we get the
sequence
$$0\to p_*\cal I_\C(6)\to p_*\oo_{\HWP_4\times\WP}(6)
\to p_*\oo_\C(6)\to R^1p_*\cal I_\C(6)\ .$$
If we assume for a moment that $R^1p_*\cal I_\C(6)$ vanishes and
that the (zeroth) direct images are locally free, this sequence reduces to
$$0\to p_*\cal I_\C(6)\to \Coh^0(\oo_{\HWP_4\times\WP}(6))_{\HWP_4}
\stackrel{\rho}{\to} p_*\oo_\C(6)\to 0\ .$$
Now take a section of $\Coh^0(\oo_{\HWP_4\times\WP}(6))_{\HWP_4}$ which is
induced
from a generic section $F$ of $\oo_{\HWP_4\times\WP}(6)$ and so
represents a Calabi-Yau hypersurface $X_F$.
If $\lbrack C\rbrack\in\HWP_4$ is a given curve, then the induced section
$\rho(F)$ of $p_*\oo_\C(6)$ restricted to the fiber
$\Coh^0(C, \oo_C(6))$ over $\lbrack C\rbrack\in\HWP_4$ is equal to the
restriction of $F$ to $C$. Hence, $\rho(F)$ vanishes exactly when $C$ is
contained in $X_F$.

Since the rank of $p_*\oo_\C(6)$ equals the dimension of $\HWP_4$,
Kleiman's Bertini theorem (cf.\ \cite{kleiman:bertini}, Remark 6) implies
that the zero scheme of the section $\rho(F)$ is finite and nonsingular;
hence the length of this scheme is equal to the number of rational quartics
on a generic Calabi-Yau hypersurface in $\WP$. This number is given by the
integral
\begin{equation}\label{integral}
\int_{\HWP_4}c_{13}(p_*\oo_\C(6))\ ,
\end{equation}
$c_{13}$ being the top Chern class.

It remains to prove the claimed facts about the direct image sheaves.

\begin{prop}
$p_*\cal I_\C(6)$ and $p_*\oo_\C(6)$ are locally free
sheaves, and $R^1p_*\cal I_\C(6)$ vanishes.
\end{prop}

\begin{pf}
We show first that $\Coh^i(\cal I_C(6)) = 0$ for all curves
$\lbrack C\rbrack\in\HWP_4$, $i\geq 1$, and that $\dim\Coh^0(\cal I_C(6))$
is constant on $\HWP_4$.

Since $\dim\Coh^i(\WP,\I_C(6))$ is an upper semicontinuous function on
$\HWP_4$, we have to show the vanishing of $\Coh^i(\WP,\I_C(6))$ only
for all degenerations with monomial ideals, and the constantness
of $\dim\Coh^0(\I_C(6))$ for those and a generic curve.

Namely, let $\lbrack C\rbrack\in\HWP_4$ and suppose
$\dim\Coh^i(\cal I_C(6)) = d$.
Take a one-dimensional torus action on $\WP$ such that the monomial  curves
are the fixed points of the action.
For all $t\in \bold C^*$, the schemes $C_t := tC$ are projectively equivalent,
hence the cohomology groups have all the same dimension $d_t = d$.
But the limit $C _0 = \lim_{t\to 0}C_t$ is a monomial curve, and by
semicontinuity,
$\dim\Coh^i(\I_{C_0}(6)) \geq \dim\Coh^i(\I_C(6)) $.

As mentioned before, all the ideals of monomial degenerations of
elliptic quartics $C'\subseteq\bold P^4$ are at least 3-regular, thus
$\Coh^i(\bold P^4,\cal I_{C'}(6)) = 0$ for all $i>0$.
By \lemref{quotients}, it follows that $\Coh^i(\WP,\I_C(6)) = 0$
for all $i>0$ and all monomial degenerations of rational quartics.
The constantness of $\dim\Coh^0(\I_C(6))$ can also be verified by
an explicit computation.

Since $\oo_\C(6)$ is flat over $\HWP_4$, $\cal I_\C(6)$
is a flat sheaf, too. By the previous result and cohomology and base change
theorems (\cite{thehartshorne}, III 12.11, 12.9) we conclude that
$R^1p_*\I_\C(6) = 0$ and that
$p_*\I_\C(6)$ and $p_*\oo_\C(6)$ are locally free.
\end{pf}

\section{TheCalculation}
\label{calc}
We will calculate the integral \eqref{integral} by Bott's formula, as
follows (cf. \cite{bott:formula} and \cite{geirsas:twcub2}; these ideas are
largely due to Geir Ellingsrud and Stein Arild Str{\o}mme):

Suppose we are given a $\bold C^*$-action on $\WP$ which induces a $\bold
C^*$-action with isolated fixed points on $\HWP_4$. This action in turn
induces an equivariant $\bold C^*$-action on the tangent bundle
$\cal T_{\HWP_4}$ and on $p_*\oo_\C(6)$. Therefore, in a fixed point
$x = \lbrack C \rbrack$ of the action, the respective fibers $\cal
T_{\HWP_4}(x)$
and $p_*\oo_\C(6)\otimes\bold C(x)$ are $\bold C^*$-representations.
As torus representations, they decompose into a direct sum of one-dimensional
representations. Let $w_1(x),\ldots,w_r(x)$ resp.\ $\tau_1(x),\ldots,\tau_r(x)$
be the corresponding weights.
\goodbreak
Then Bott's formula says in our context:
\begin{equation}\label{bott}
\int_{\HWP_4} c_{13}(p_*\oo_\C(6)) =
\sum_{x\in\HWP_4^{\bold C^*}}{\tau_1(x)\cdot\ldots\cdot\tau_r(x)\over
w_1(x)\cdot\ldots\cdot w_r(x)}.
\end{equation}
Let $T\subseteq \GL(5)$ be a maximal torus which acts diagonally
on $\WP$ w.r.t.\ the coordinates $x_0,\ldots,x_4$ of $\bold P^4$.
There are characters $\lambda_0,\ldots,\lambda_4$ on $T$ such that for any
$t\in T$,
we have $t\cdot x_i = \lambda_i(t)\cdot x_i$,
and these characters generate the representation ring of $T$, i.e.,
if $W$ is a finite representation of $T$, we can write {\it cum granum salis:}
\goodbreak
$$W = \sum a_{p_0\ldots p_4}\lambda_0^{p_0}\cdot\ldots\cdot\lambda_4^{p_4}\ .$$
The action of $T$ descends to an action on $\WP = \bold P^4/\Gamma$.

In the following, we will compute the torus representations of the induced
$T$-action on $\HWP_4$ in the fibers of $p_*\oo_\C(6)$ and
$\cal T_{\HWP_4}$ in all fixed points. It is easy to see that a point
$x\in\HWP_4$ is fixed exactly when the graded ideal of the corresponding
curve is generated by monomials.

Then we choose a one-parameter subgroup $\bold C^*\subseteq T$ with no
non-trivial $\bold C^*$-weight in the tangent space of any fixed point.
Such a one-parameter subgroup is given by a point $(w_0,\ldots,w_4)$
in the weight lattice $\Hom(\bold C^*,T)\cong\bold Z^5$; the
corresponding characters on $\bold C^*$ are given by $\lambda_i(t)=t^{w_i}$.
If $\lambda_0^{p_0}\cdot\ldots\cdot\lambda_4^{p_4}$ is the character
of the $\bold C^*$-representation on an invariant one-dimensional
subspace of the tangent space in a fixed point, the corresponding
weight is given by
\begin{equation}\label{weight}
w=p_0w_0+\ldots+p_4w_4.
\end{equation}
All these weights are nonzero if the weight vector $(w_0,\ldots,w_4)$
is chosen to avoid simultaneously all the (finitely many) hyperplanes
in the weight lattice defined by the linear forms \eqref{weight}.
Such a choice is clearly possible.
(In the concrete calculation of the integral \eqref{bott} by the
{\sc Maple}-program listed in the appendix, we
try randomly chosen weights; if none of the denominators in the
summands of \eqref{bott} is zero --- which would result in a
``division by zero" error --- the choice is valid.)

Our choice of the weights guarantees that the $\bold C^*$-action on $\HWP_4$
has isolated fixed points (in fact, the same fixed points as the action
of $T$), hence we will be able to apply Bott's formula.\\

We will first calculate the fixed points and tangent space representations
for $\HWP_3$ and afterwards use the fact that $\HWP_4$ is a locally trivial
fibration over $\check{\bold P}^3$ with fiber~$\HWP_3$.

This is being done by calculating the fixed points of the $T$-action
and the $T$-representation on the tangent spaces in each successive step
of the blow up, i.e.\ on $\tilde G$, on $\tilde G'$, and finally on $\HWP_3$
(Note that there are induced $T$-actions on those spaces and that
fixed points lie over fixed points).

To compute the data of a blow up, say, of $X$ with center $B$, we have to get
hold of
the fiber  $\cal N(x)$ of the normal bundle $\cal N=\cal N_{B/X}$ to the
subvariety $B$ to blow up in each fixed point $x\in B$. We can achieve this by
computing $T$-semiinvariant base vectors of $\cal N(x)$. Every such vector
$\xi$
gives rise to a fixed point $x_\xi$ in the proper transform $B'$ of $B$:
$x_\xi$ is the inverse image of $x$ in the proper transform of a curve
tangent to $\xi$ in $x$.

The tangent space at $x_\xi$ in the blow up $X' = Bl_B X$ is then given
as
\begin{equation}\label{tangentofblowup}
\cal T_{X'}(x_\xi)
\cong L_\xi\oplus\cal T_B(x)\oplus\cal T_{\bold P(\cal N^\vee(x))}(x_\xi)\ ,
\end{equation}
where $L_\xi$ is the span of $\xi$ in the normal space $\cal N(x)$.
The isomorphism is equivariant.
(Recall that the exceptional divisor of the blow up with center $B$ is
isomorphic to the projective bundle $\bold P(\cal N^\vee_{B/X})$, with
normal space in $x_\xi$ isomorphic to $L_\xi$;
cf.\ for instance \cite{fulton:intersect}, B.~6.) \\

Let us now look at the concrete calculations.
We denote by $\tilde B'$ (resp.\ $\tilde D'$) the proper transform of
$\tilde B$ in $\tilde G$ (resp.\ of $\tilde D$ in $\HWP_3$).

\begin{prop}
There are 126 fixed points in $\HWP_3$, and they are projectively
equivalent to one of the following 25 types listed below. (To each
fixed point type, we give the permutations of the variables which generate
the remaining fixed points of that type.)
\begin{equation*}
\def\ulspace{\vrule height12pt depth0pt width0pt}
\def\llspace{\vrule height0pt depth6pt width0pt}
\begin{tabular}{|c|c|c|c|}
\hline\ulspace{\bf lies} & & & {\bf number of} \\
{\bf in}\llspace & \raisebox{1.5ex}[0cm][0cm]{\bf fixed point type} &
\raisebox{1.5ex}[0cm][0cm]{\bf permutations} & {\bf fixed points}    \\
\hline
& \ulspace $(x_0^2,x_1^2)$ & cyclic & 3\\
\raisebox{-1.5ex}[0cm][0cm]{$\tilde G-\tilde B$}
& \ulspace $(x_0^2,x_1x_2)$ & permutations & 3\\
& \ulspace $(x_1^2,x_2^2)$ & of & 3\\
& \ulspace\llspace $(x_1x_2,x_3^2)$ & $x_1,x_2,x_3$  & 3\\
\hline
&\ulspace $(x_1x_2,x_1x_3,f),$ & $x_1 \leftrightarrow x_2$ &
\raisebox{-1.5ex}[0cm][0cm]{6 $\times$ 4}\\
\raisebox{-1.6ex}[0cm][0cm]{$\tilde B'-\tilde D$}
& \llspace $f\in \{x_0^2x_2, x_0^2x_3, x_2^3, x_2^2x_3, x_2x_3^2, x_3^3\}$ &
       and $x_1 \leftrightarrow x_3$ &\\
& \ulspace $(x_1^2,x_1x_2,g),$ & all permutations &
\raisebox{-1.5ex}[0cm][0cm]{3 $\times$ 6}\\
& \llspace $g\in\{x_2^3, x_22x_3, x_2x_3^2, x_0^2x_2, x_1x_3^2, x_1x_0^2\}$
       &  of $x_1,x_2,x_3$ &\\
\hline
& \ulspace  $(x_1^2, x_1x_2, x_1x_3^2,f),$ &
\raisebox{-3pt}[0cm][0cm]{all} &
\raisebox{-1.5ex}[0cm][0cm]{6 $\times$ 6} \\
\raisebox{-1.5ex}[0cm][0cm]{$\tilde D'$}
& \llspace $f\in\{x_2x_3^3, x_2^2x_3^2, x_3^4, x_2^4, x_2^3x_3, x_0^2x_2^2\}$ &
       \raisebox{-2pt}[0cm][0cm]{permutations} &\\
& \ulspace $(x_1^2, x_1x_2, x_1x_0^2,g),$ & \raisebox{3pt}[0cm][0cm]{of}
 & \raisebox{-1.5ex}[0cm][0cm]{6 $\times$ 6} \\
& \llspace $g\in \{x_0^2x_2x_3, x_2^2x_3^2, x_0^4, x_2^4, x_2^3x_3,
       x_0^2x_2^2\}$ & \raisebox{4pt}[0cm][0cm]{$x_1,x_2,x_3$} &\\
\hline
\end{tabular}
\end{equation*}
The tangent space in a fixpoint $x$ in $\tilde G$ is given by
$\cal T_{\tilde G}(x) =
\sum \lambda_\gamma\lambda_\delta\lambda^{-1}_\alpha\lambda^{-1}_\beta,$
where $\alpha$, $\beta$ ($\alpha\leq\beta$) run over the pairs of indices
of the monomials in $I_x$, whereas $\gamma$, $\delta$ ($\gamma\leq\delta$)
run over all except these indices.
The tangent space in each fixed point $x_\xi$ lying over a fixed point $x$
in $\tilde B$ resp.\ $\tilde D$ is calculated by formula
\eqref{tangentofblowup},
where the tangent spaces and normal spaces to $\tilde B$ resp.\ $\tilde D$
are given by:
\begin{equation*}
\def\ulspace{\vrule height12pt depth0pt width0pt}
\def\llspace{\vrule height0pt depth6pt width0pt}
\begin{tabular}{|c|rl|}
\hline\ulspace\llspace $I_{x_\xi}$ & {\bf tangent\!\!\!\!\!} & {\ \bf and
normal spaces} \\
\hline\ulspace\llspace &
$\cal T_{\tilde B}(x)\quad =$ &
$\lambda_1\lambda_2^{-1} +
\lambda_1\lambda_3^{-1} +
\lambda_2\lambda_1^{-1} +
\lambda_3\lambda_1^{-1}$ \\
 $(x_1x_2,x_1x_3,f)$ & $\cal N_{\tilde B/\tilde G}(x)\quad =$ &
$\lambda_0^2\lambda_1^{-1}\lambda_2^{-1} +
\lambda_0^2\lambda_1^{-1}\lambda_3^{-1} +
\lambda_2\lambda_1^{-1} +{}$\\
\llspace & &$\lambda_3\lambda_1^{-1} +
\lambda_3^2\lambda_1^{-1}\lambda_2^{-1} +
\lambda_2^2\lambda_1^{-1}\lambda_3^{-1} $ \\
\ulspace\llspace &
$\cal T_{\tilde B}(x)\quad =$ &
$\lambda_3\lambda_1^{-1} +
\lambda_3\lambda_2^{-1} +
\lambda_2\lambda_1^{-1} +
\lambda_3\lambda_1^{-1}$ \\
$(x_1^2,x_1x_2,g)$ & $\cal N_{\tilde B/\tilde G}(x)\quad =$ &
$\lambda_0^2\lambda_1^{-1}\lambda_2^{-1} +
\lambda_0^2\lambda_1^{-2} +
\lambda_3^2\lambda_1^{-1}\lambda_2^{-1} +{}$\\
\llspace& &$\lambda_3^2\lambda_1^{-2} +
\lambda_2\lambda_3\lambda_1^{-2} +
\lambda_2^2\lambda_1^{-2} $ \\
\hline
\ulspace\llspace &
$\cal T_{\tilde D}(x)\quad =$ &
$\lambda_3\lambda_1^{-1} +
\lambda_2\lambda_1^{-1} +
\lambda_3\lambda_2^{-1} +
\lambda_0^2\lambda_3^{-2}$ \\
$(x_1^2, x_1x_2, x_1x_3^2,f)$ & $\cal N_{\tilde D/\tilde G'}(x)\quad =$ &
$\lambda_3\lambda_1^{-1} +
\lambda_3^2\lambda_1^{-1}\lambda_2^{-1} +
\lambda_2^3\lambda_1^{-1}\lambda_3^{-2} +{}$\\
\llspace & &$\lambda_2^2\lambda_1^{-1}\lambda_3^{-1} +
\lambda_2\lambda_1^{-1} +
\lambda_0^2\lambda_2\lambda_1^{-1}\lambda_3^{-2} $ \\
\ulspace\llspace &
$\cal T_{\tilde D}(x)\quad =$ &
$\lambda_3\lambda_1^{-1} +
\lambda_2\lambda_1^{-1} +
\lambda_3\lambda_2^{-1} +
\lambda_3^2\lambda_0^{-2}$ \\
$(x_1^2, x_1x_2, x_1x_0^2,g)$ & $\cal N_{\tilde D/\tilde G'}(x)\quad =$ &
$\lambda_3\lambda_1^{-1} +
\lambda_0^2\lambda_1^{-1}\lambda_2^{-1} +
\lambda_2^3\lambda_1^{-1}\lambda_0^{-2} +{}$\\
\llspace & & $\lambda_2^2\lambda_3\lambda_0^{-2}\lambda_1^{-1} +
\lambda_2\lambda_1^{-1} +
\lambda_2\lambda_3^2\lambda_0^{-2}\lambda_1^{-1} $ \\
\hline
\end{tabular}
\end{equation*}

\end{prop}
\begin{pf}
The fixed points $x$ in the grassmannian $\tilde G$ are readily determined,
their ideals have the form $I_x = (x_ix_j, x_kx_l)$ with the obvious
restrictions on the indices.
Let $V_2 := \Coh^0(\oo_{\WP^3}(2))$,
$V_{I_x} := \bold C\cdot x_ix_j\oplus\bold C\cdot x_kx_l$.
The tangent space in a fixed point $x$ is given by
$$\cal T_{\tilde G}(x) = \Hom(V_{I_x}, V_2/V_{I_x}) =
\sum \lambda^{-1}_\alpha\lambda^{-1}_\beta\lambda_\gamma\lambda_\delta\ ,$$
the indices being as specified in the proposition.\\

{\it First blow up.}
The subvariety $\tilde B\subseteq\tilde G$ to blow up consists of pencils
with a fixed component, and the fixed points in $\tilde B$ are of
type $(x_1x_2, x_1x_3)$ and $(x_1^2, x_1x_2)$.

Consider the fixed point $x$ with ideal $I_x = (x_1x_2, x_1x_3)$.
The tangent space to $\tilde G$ in $x$ is given by
\begin{eqnarray}\label{tangentfistbl}
\cal T_{\tilde G}(x) &=
&\lambda_1^{-1}\lambda_2^{-1}\lambda_0^2 +
\lambda_1^{-1}\lambda_3^{-1}\lambda_0^2 +
\lambda_2^{-1}\lambda_1 +
\lambda_3^{-1}\lambda_1 +
2\lambda_1^{-1}\lambda_2 +\\
&&2\lambda_1^{-1}\lambda_3 +
\lambda_1^{-1}\lambda_2^{-1}\lambda_3^2 +
\lambda_1^{-1}\lambda_3^{-1}\lambda_2^2\ .\nonumber
\end{eqnarray}
First, we will determine a semiinvariant basis for the fiber
$\cal N_{\tilde B/\tilde G}(x)$ of the normal bundle of $\tilde B$
in $\tilde G$.
Let $\xi\in \cal T_{\tilde G}(x)$ be a semiinvariant tangent vector in $x$,
given as $\xi = {\xi_1\choose\xi_2}$ w.r.t.\ the basis
$\{x_1x_2,x_1x_3\}$ of $V_{I_x}$
(i.e., $\xi(ax_1x_2,bx_1x_3) = (a\xi_1x_1x_2,b\xi_2x_1x_3)$).
Since $\xi$ is semiinvariant, $\xi_1$ and $\xi_2$ are scalar multiples of a
common Laurent monomial $\mu_\xi$ in $x_0,\ldots,x_3$ of degree 0.
Furthermore, the torus representation on the subspace spanned by $\xi$ is
obtained by formally substituting $\lambda_i$ for $x_i$, $i=0,\ldots,3$
in $\mu_\xi$. Clearly this monomial in the $\lambda_i$ has to be one of the
summands in \eqref{tangentfistbl}.

Now $I_\xi(t) = (x_1x_2+t\xi_1x_1x_2, x_1x_3+t\xi_2x_1x_3)$ is a curve
through $x$ with tangent direction $\xi$ in $x$.
We see that the semiinvariant tangent vectors
$$\xi =
{x_1x_2^{-1}\choose 0},
{0\choose x_1x_3^{-1}},
{x_2x_1^{-1}\choose x_2x_1^{-1}},
{x_3x_1^{-1}\choose x_3x_1^{-1}}
$$
are tangent to $\tilde B$; the curves given by $I_\xi(t)$ are even contained in
$\tilde B$.
The vectors
$$
{x_2x_1^{-1}\choose -x_2x_1^{-1}},
{x_3x_1^{-1}\choose -x_3x_1^{-1}},
{x_0^2x_1^{-1}x_2^{-1}\choose 0},
{x_3^2x_1^{-1}x_2^{-1}\choose 0},
{0 \choose x_0^2x_1^{-1}x_3^{-1}},
{0 \choose x_2^2x_1^{-1}x_3^{-1}}
$$
complete the previous ones to a semiinvariant basis of $\cal T_{\tilde G}(x)$,
thus they represent (modulo $\cal T_{\tilde B}(x)$) a semiinvariant
basis of $\cal N_{\tilde B/\tilde G}(x)$.

In order to compute the fixed points of $\tilde G'$ lying over $x$,
we consider the curves $I_\xi(t)$ for $\xi$ in that basis.
Each of them defines a flat family of curves over $\bold A^1-\{0\}$.
We can extend this family in a unique way to a flat family over $\bold A^1$.
That flat family induces a map $\bold A^1\to\HWP_3$ such that the
image $x_\xi'$ of $0\in\bold A^1$ maps onto a fixed point $x_\xi$
of $\tilde G'$ under the blow up map $\HWP_3\to\tilde G'$.
The ideal $I_{x_\xi}$ corresponding to $x_\xi$ is the subideal generated in
degree 3 of the ideal corresponding to $x_\xi'$ (cf.\ \eqref{H3toG'}).

To actually compute the ideal $I_{x_\xi}$, we use a flattening algorithm
described by Bayer and Mumford (\cite{baymum:compute}, Ch.1):
\\
\begin{prop}
Let $I(t) = (m_1(t),\ldots,m_r(t))$ be the ideal of a family of schemes over
$\bold A^1$ and suppose that all $m_i := m_i(0)$ are monomials.
Let $I=(m_1,\ldots,m_r)$ be the ideal of the central fiber.
Consider the following algorithm:
\begin{enumerate}
\item
Take a minimal syzygy of two generators $m_i$, $m_j$ of the ideal $I$,
i.e., a relation
$$h_im_i+h_jm_j=0,$$
with $h_i$ and $h_j$ coprime monomials,
which does not lift to a syzygy in $I_\xi(t)$.
This means that we get a relation
$$h_im_i(t)+h_jm_j(t) = t^sg$$
with $g\neq 0$ and $t \not\:\mid g$.
Add the polynomial $g$ to the generators of $I(t)$
to get new ideals $I'(t)$ and $I'=I'(0)$.
\item
Repeat this process finitely often until all syzygies lift to syzygies
in $I'(t)$.
\end{enumerate}
Then the resulting ideal $I'(t)$ is flat in $t=0$, and the ideal
of the fiber is equal to~$I'$.
\end{prop}
In our case, all $g$'s are monomials and $s$ is equal to 1 in the first step.
Furthermore, we are only interested in generators of degree 3, so we can stop
the iteration when we have added all monomials of degree 3, and that is the
case
already after the first step, as can easily be seen.

In concrete terms, the only syzygy $x_3\cdot(x_1x_2)-x_2\cdot(x_1x_3)=0$
lifts to $I_\xi(t)$ as
\begin{eqnarray*}
x_3\cdot(x_1x_2+t\xi_1x_1x_2)-x_2\cdot(x_1x_3+t\xi_2x_1x_3)
&= &t\cdot(x_3\xi_1x_1x_2-x_2\xi_2x_1x_3)\\
&=: &t\cdot f_\xi\ .
\end{eqnarray*}

All together,
the fixed points in the blow up $\tilde G'$ lying over $x$ are:
$$I_\xi = (x_1x_2,x_1x_3,f_\xi), \quad
f_\xi\in \{x_0^2x_2, x_0^2x_3, x_2^3, x_2^2x_3, x_2x_3^2, x_3^3\}\ ,$$
plus the ideals obtained by the permutations $x_1 \leftrightarrow x_2$
and $x_1 \leftrightarrow x_3$.
The  calculation for points $x$ of type $I_x = (x_1^2, x_1x_2)$
is exactly analogous.\\

{\it Second blow up.}
In $\tilde G'$, we have to blow up the non-flat locus
$\tilde D\subseteq \tilde G'$. The fixed points $x$ of $\tilde G'$ contained
in $\tilde D$ are of type $I_x = (x_1^2, x_1x_2, x_1x_3^2)$ and
$I_x = (x_1^2, x_1x_2, x_1x_0^2)$.

The tangent space to $\tilde G'$ in the first fixed point is given by
(cf.\ \eqref{tangentofblowup}):
\begin{eqnarray}\label{tangentsecondbl}
\cal T_{\tilde G}(x) &=
&2\lambda_1^{-1}\lambda_3 +
2\lambda_1^{-1}\lambda_2 +
\lambda_2^{-1}\lambda_3 +
\lambda_3^{-2}\lambda_0^2 +
\lambda_1^{-1}\lambda_2^{-1}\lambda_3^2 +\\
&&\lambda_1^{-1}\lambda_3^{-2}\lambda_2^3 +
\lambda_1^{-1}\lambda_3^{-1}\lambda_2^2 +
\lambda_1^{-1}\lambda_3^{-2}\lambda_2\lambda_0^2\ .\nonumber
\end{eqnarray}
We will determine (the torus representation of) the normal space $\cal N(x)$ to
$\tilde D$ in $x$.

%Since $\tilde G'$ is naturally embedded in
%$G_3 := \Grass_6\bigl(\Coh^0(\oo_{\bold P(2,1^3)}(3))\bigr)$,
%a tangent vector $\xi\in\cal T_{\tilde G'}(x)$ can be viewed as an element of
%$\cal T_{G_3}(x) = \Hom(V_{I_x}, V_3/V_{I_x})$, where
%$V_3 = \Coh^0(\oo_{\WP^3}(3))$ and
%\begin{eqnarray*}
%V_{I_x} &= &span(\hbox{monomials of deg. 3 in } I_x)\\
%&= &span(x_1^3,x_1^2x_2,x_1^2x_3,x_1x_2^2,x_1x_2x_3,x_1x_3^2)\ .
%\end{eqnarray*}

A curve in $\tilde G'$ through $x$ tangent to $\xi\in\cal T_{\tilde G'}(x)$
is given  modulo $t^2$ as
$$I_\xi(t) = (x_1^2+t\xi_1x_1^2, x_1x_2+t\xi_2x_1x_2,
x_1x_3^2+t\xi_3x_1x_3^2)\ ,$$
where $\xi = (\xi_1, \ldots,\xi_3)^t$, and if $\xi$ is semiinvariant, then
the $\xi_i$ are scalar multiples of a common Laurent monomial $\mu_\xi$ of
degree zero.

Again, by lifting a syzygy relation from $I_x$ to $I_\xi(t)$ (this can be
done modulo $t^2$, too),
we calculate the monomial $f_\xi$ we have to add to $I_x$ in order to
get the fixed point in the blow up $\HWP_3$ corresponding to $\xi$.
Since the syzygies of all pairs of monomials in $I_x$ generate the whole
syzygy module, we only need to consider such pairs.

Suppose that lifting the syzygy of the pair $(m_1, m_2)$ results in the
right monomial $f_\xi$; then $f_\xi$ is equal to
$lcm(m_1,m_2)\cdot\mu_\xi$.
This monomial is supposed to have degree four. The two syzygies which can
possibly yield a monomial of degree four are those between $x_2$ and $x_1x_3^2$
and between $x_1x_2$ and $x_1x_3^2$. Consider the first pair. We have
$f_\xi = x_1^2x_3^2\cdot\mu_\xi$. But from \eqref{tangentsecondbl} it is
apparent that the Laurent monomials $\mu_\xi$ contain $x_1$ to the power -1 or
higher, that means that $f_\xi$ contains $x_1$ as a factor. Hence the
resulting ideal $I_\xi=(x_1^2,x_1x_2,x_1x_3^2,f_\xi)$ is contained in the
ideal $(x_1)$ and certainly doesn't correspond to a flat degeneration.

So we only have to consider the pair $x_1x_2$, $x_1x_3^2$.
We claim that the monomials $\mu_\xi$ corresponding to a semiinvariant basis
of $\cal N(x)$ are
\begin{eqnarray}\label{mus}
\{x_1^{-1}x_3,\  x_1^{-1}x_2,\  x_1^{-1}x_2^{-1}x_3^2,\  x_1^{-1}x_3^{-2}x_2^3,
\ x_1^{-1}x_3^{-1}x_2^2,\  x_1^{-1}x_3^{-2}x_2x_0^2\}\ .
\end{eqnarray}

The two other possibilities for $\mu_\xi$, namely $x_2^{-1}x_3$
and $x_3^{-2}x_0^2$, are excluded because they again lead to ideals contained
in the ideal $(x_1)$.
On the other hand, $\cal N(x)$ has dimension 6; together this proves the claim.

So, the fixed points in the blow up $\HWP_3$ lying over
$I_x = (x_1^2, x_1x_2, x_1x_3^2)$
are
$$I_\xi = (x_1^2, x_1x_2, x_1x_3^2,f_\xi), \quad
f_\xi\in \{x_2x_3^3, x_2^2x_3^2, x_3^4, x_2^4, x_2^3x_3, x_0^2x_2^2\}\ .$$
The torus representation of $\cal N(x)$ is gotten by formally summing up the
monomials in \eqref{mus} and substituting $\lambda$ for $x$.
To get all fixed points, we have to add the fixed points generated by
all permutations of the variables $x_1$, $x_2$ and $x_3$.
The calculation for points $x$ of type
$I_x = (x_1^2, x_1x_2, x_1x_0^2)$ is again analogous.

Now we have determined all fixed points of the torus action on $\HWP_3$:
they are the fixed points of all stages of the blow up process minus the
fixed points lying in the blow up loci.
\end{pf}

{\it Fixed points and tangent spaces in $\HWP_4$.}
A fixed curve in $\WP$ spans a $T$-invariant hyperplane
($\cong \bold P(2,1^3)$) in $\WP$, which is given by an equation
$x_i = 0$, $i\in\{1,\ldots,4\}$.
Thus, if for instance $i=4$, we get the ideals of all fixed curves
lying in $\{x_4=0\}$ by adjoining the monomial $x_4$ to the ideals
previously determined. The ideals of the remaining fixed points are
obtained by cyclically permuting the variables $x_1,\ldots,x_4$ in these
ideals. Thus, $\HWP_4$ containes 504 fixed points all together.

According to \propref{fiberingofhilb}, $\HWP_4$ is fibered over
$\check{\bold P}^3$ as
\begin{eqnarray*}
\pi: \HWP_4 &\to &\check{\bold P}^3 \\
\lbrack C \rbrack & \mapsto & \lbrack\hbox{ projection of }
span(C) \hbox{ to } \bold P^3\subseteq\WP\rbrack\ .
\end{eqnarray*}
The tangent space of $\HWP_4$ in a fixed point $x=\lbrack C \rbrack$
therefore decomposes equivariantly as a direct sum
$$\cal T_{\HWP_3}(x)\oplus\cal T_{\check{\bold P}^3}(\pi(x)) .$$
Let $C$ be an invariant curve spanning the singular hyperplane $\{x_i=0\}$
and let $V_i = \bold C\cdot x_i$,
$V = \bold C\cdot x_1\oplus\ldots\oplus\bold C\cdot x_4$.
Then the torus representation of the tangent space to $\check{\bold P}^3$
is equal to
$$\cal T_{\check{\bold P}^3}(\pi(x)) = \Hom(V_i, V/V_i)
= \sum_{j\neq i}\lambda_j\lambda_i^{-1}\ $$

Finally, the torus representations on the fibers $\Coh^0(\oo_C(6))$
of $p_*\oo_\C(6)$ are easily determined:
$\Coh^0(\oo_C(6))$ is spanned by those monomials of degree 6 that don't lie
in the ideal $\cal I_C$.

%\newpage
\section{Appendix: Program Listing}
\begin{verbatim}
# Calculation of the number of rational quartics on a general
# Calabi-Yau hypersurface in P(2,1,1,1,1).

### some setup and utility routines

# is_invariant decides whether a monomial is invariant under
# var -> -var, i.e. whether the variable var occurs in even power

is_invariant := proc(mon, var)
   mon = subs(var=-var,mon)
end:

# invariant_subspace returns the invariant subspace of a vector space
# (w.r.t. var -> -var)

invariant_subspace := proc(J, var) local i, inv;
   inv := 0:
   for i from 1 to nops(J) do
      if is_invariant(op(i,J),var) then
         inv := inv + op(i,J)
      fi
   od:
   inv:
end:

# setuptorus(n,m) sets up the data of an (n+1)-dimensional torus
# with characters x0,...xn. The representation ring is identified
# with the Laurent series ring Z[xi,1/xi, i=0..n]. V[i] denotes the
# canonical representation of H^0(o_P(2,1^n)(i)).

setuptorus := proc(N,M) local genf,i,k;
   genf := expand(convert(series(
         1/product('1-x.i*t',i=0..N),t,M+1),polynom)):
   for i from 0 to M do
      V[i] := invariant_subspace(sort(coeff(genf,t,i)),x0)
   od;
   i := 'i':
   dualsubstring := seq(x.i=x.i^(-1),i=0..N):
   variables := [seq(x.i,i=0..N), seq(x.i^(-1),i=0..N)]:
end:

dualrep := proc(C) sort(expand(subs(dualsubstring,C))) end:

simpleweights := proc(p) local mon;
   expand(p); if p=0 then RETURN(0) fi;
   coeffs(",variables,'mon');
   convert([mon],`+`);
end:

# twist(J,k) computes the representation of the degree k part of
# the ideal J:

twist := proc(J,k) local i,j,res,mon,deg;
   res := 0:
   for j from 1 to nops(J) do
      mon := op(j,J):
      if k >= degree(mon) then
         res := res + V[k-degree(mon)]*mon:
      fi:
   od:
   simpleweights(res):
end:

# If C* -> T is the one-parameter subgroup given by weights
# w0,..,w4, the character xi restricts to t^(wi), i=0..4.
# In order to calculate the enumerator in Bott's formula, we have
# to calculate the product of all of them.

prodwts := proc(f) local t,cof,mon,res,i;
   cof := [coeffs(f,variables,mon)];
   mon := subs(x0=t^w0,x1=t^w1,x2=t^w2,x3=t^w3,[mon]);
   res := 1;
   for i from 1 to nops(mon) do
      res := res*subs(t=1,diff(mon[i],t))^cof[i]
   od;
   res
end:


### here begins the actual calculation

# First, we calculate all data for H_3. After that, we use the
# fiber structure of H_4 with fiber H_3 to calculate the torus
# action on H_4.

setuptorus(3,6);

tw := 6:    # we are going to calculate o(6)
fp := 0:      # fp counts the number of fixpoints

# ideals and tangent spaces of the fixed points in the grassmannian

for i0 from 0 to 3 do
 for i1 from i0 to 3 do
  for i2 from i0 to 3 do
   for i3 from i2 to 3 do
      if i0 <> i2 and i0 <> i3 and i1 <> i2 and i1 <> i3 then
         idd := x.i0*x.i1 + x.i2*x.i3:     # ideal of fixed point
         if is_invariant(idd,x0) then      # we have to filter
            fp := fp+1:                    #  out the invariant
            id[fp] := idd:                 #  ideals
            T[fp] := expand((V[2]-id[fp])  # tangent space of
                     *dualrep(id[fp])):    #  the grassmannian
            o[fp] :=                       # torus representation of
               V[tw]-twist(id[fp],tw):     #  the fiber of p_*O_C(6)
         fi
      fi
   od
  od
 od
od:

# first blow up

idd := x1*x2+x1*x3:                        # ideal of the first type

TB := expand((x2+x3)/x1 + x1/x2 + x1/x3):  # tangent space to the
                                           #  subvariety B to blow up
N := expand((V[2]-idd)*dualrep(idd)) - TB: # normal space to B

for i from 1 to nops(N) do                # normal direction to the
   TL := op(i,N):                         #  exceptional divisor
   TPN := expand((N-TL)/TL):   # tangent space to the fiber
   T_[i] := TB + TL + TPN:                # total tangent space
   T_[i] := unapply(T_[i], x1,x2,x3):
   id_[i] := idd+x1*x2*x3*TL:             # ideal of the fixed point
   o_[i] := V[tw]-twist(id_[i],tw):       # torus representation of
   o_[i] := unapply(o_[i], x1,x2,x3):     #  the fiber of p_*O_C(6)
   id_[i] := unapply(id_[i], x1,x2,x3):   #  in a fixed point
od:


# permutations of the variables that generate the remaining
# fixed points
permut1 := [x1,x2,x3], [x2,x1,x3], [x3,x2,x1]:

for j from 1 to 3 do                     # the remaining fixed points
   for i from 1 to nops(N) do
      fp := fp+1:
      T[fp] := T_[i](op(permut1[j])):
      o[fp] := o_[i](op(permut1[j])):
      id[fp] := id_[i](op(permut1[j])):
   od:
od:

idd := x1^2+x1*x2:                          # ideal of the second type
TB := expand(x3*dualrep(x1+x2) + (x2+x3)*dualrep(x1)):
N := expand((V[2]-idd)*dualrep(idd)) - TB:

for i from 1 to nops(N) do
   TL := op(i,N):
   TPN := expand((N-TL)/TL):
   T_[i] := TB + TL + TPN:
   T_[i] := unapply(T_[i], x1,x2,x3):
   id_[i] := idd+x1^2*x2*TL:
   if divide(id_[i],x1) then o_[i] := 0 else # the ideals whose
      o_[i] := V[tw]-twist(id_[i],tw) fi:    # generators have a
   o_[i] := unapply(o_[i], x1,x2,x3):        # common factor are lying
   id_[i] := unapply(id_[i], x1,x2,x3):      # in the next blow up
od:                                          # locus

# permutations of the variables that generate the remaining
# fixed points

permut2 := [x1, x2, x3], [x2, x1, x3], [x1, x3, x2],
           [x3, x2, x1], [x3, x1, x2], [x2, x3, x1]:

for j from 1 to 6 do
   for i from 1 to nops(N) do
      fp := fp+1:
      T[fp] := T_[i](op(permut2[j])):
      o[fp] := o_[i](op(permut2[j])):
      id[fp] := id_[i](op(permut2[j])):
      if o[fp] = 0 then fp := fp-1 fi:
   od:
od:

# second blow up

idd := x1^2+x1*x2:
TL_ := [x3^2/(x1*x2), x0^2/(x1*x2)]:

for k from 1 to nops(TL_) do
   TD := expand(x3/x1 + x3/x2 + x2/x1
         + ((x3^2+x0^2)/(x1*x2)-TL_[k])/TL_[k]):
   N := expand(x3/x1 + TL_[k]
         + (x2^2+x2*x3+x3^2+x0^2)/(x1^2*TL_[k])):
   for i from 1 to nops(N) do
         # normal direction to the exceptional divisor:
      TL := op(i,N):
         # tangent space to the fiber:
      TPN := expand((N-TL)/TL):
         # total tangent space:
      T_[i] := TD + TL + TPN:
      T_[i] := unapply(T_[i], x1,x2,x3):
         # ideal of the fixed point:
      id_[i] := idd+x1^2*x2*TL_[k]+x1^2*x2^2*TL_[k]*TL:
         # torus representation of the fiber of p_*O_C(6)
         # in a fixed point
      o_[i] := V[tw]-twist(id_[i],tw):
      o_[i] := unapply(o_[i], x1,x2,x3):
      id_[i] := unapply(id_[i], x1,x2,x3):
   od:

   for j from 1 to 6 do
      for i from 1 to nops(N) do
         fp := fp+1:
         T[fp] := T_[i](op(permut2[j])):
         o[fp] := o_[i](op(permut2[j])):
         id[fp] := id_[i](op(permut2[j])):
         if o[fp] = 0 then fp := fp-1 fi:
      od:
   od:
od:

lprint(`number of fixpoints =`, fp):  # Euler characteristic
                                      # of H_3 = 126

# weight vector

# the weights are experimentally chosen in such a way that the C*
# representation has no trivial weight in any tangent space. That
# ensures that the C* action has isolated fixed points.

WG[1] := 4: WG[2] := 17: WG[3] := 55: WG[4] := 160: WG[0] := 267:

for i from 1 to fp do               # calculate the summands of Bott's
   num := prodwts(o[i]):            #  formula as a function of the
   den := prodwts(T[i]):            #  weights, still leaving out the
   fra[i] := unapply(num/den,         #  contribution of the base P^3
      w0,w1,w2,w3):                   #  of the fibration H_4 -> P^3
od:

perm_h := [0,2,3,4], [0,3,4,1],       # permutation of the hyperplanes
          [0,4,1,2], [0,1,2,3]:      #  in P(2,1^4) containing
                                    #  the fixed curves

result := 0:                        # summing up variable
for k from 1 to 4 do
   h := perm_h[k]:
   TG := (WG[h[2]]-WG[k])            # tangent space of P^3
      *(WG[h[3]]-WG[k])*(WG[h[4]]-WG[k]):
   for i from 1 to fp do            # summing up over all fixed points
      result := result+
         fra[i](WG[h[1]],WG[h[2]],WG[h[3]],WG[h[4]])/TG:
   od:
od:

lprint(`number of curves =`, result): # 6028452


\end{verbatim}

\bibliographystyle{amsplain}
\ifx\undefined\bysame
\newcommand{\bysame}{\leavevmode\hbox to3em{\hrulefill}\,}
\fi

\end{document}